\documentclass[a4paper]{amsart}
\usepackage[english]{babel}
\usepackage{amssymb}
\usepackage{graphicx}
\usepackage{booktabs}
\usepackage{longtable}
\usepackage[latin1]{inputenc}
\usepackage{url}
\usepackage[T1]{fontenc}

\newcommand{\Z}{\mathbb{Z}}

\newcommand{\Q}{\mathbb{Q}}

\theoremstyle{definition}

\theoremstyle{remark}

\DeclareMathOperator{\rank}{rank}

\begin{document}

\title[High rank  elliptic curves over quadratic fields]{High rank elliptic curves
with prescribed torsion group over quadratic fields}
\author[J. Aguirre]{Juli\'an Aguirre}
\address{
Departamento de Matem\'aticas\\
Universidad del Pa\'{\i}s Vasco\\
Aptdo. 644, 48080 Bilbao, Spain
}
\email[J. Aguirre]{julian.aguirre@ehu.es}
\author[A. Dujella]{Andrej Dujella}
\address{Department of Mathematics\\University of Zagreb\\Bijeni\v{c}ka cesta 30, 10000 Zagreb, Croatia}
\email[A. Dujella]{duje@math.hr}
\author[M. Juki\'c Bokun]{Mirela Juki\'c Bokun}
\address{Department of Mathematics\\University of Osijek\\Trg Ljudevita Gaja 6, 31000 Osijek,  Croatia}
\email[M. Juki\'c Bokun]{mirela@mathos.hr}
\author[J. C. Peral]{Juan Carlos Peral}
\address{Departamento de Matem\'aticas\\Universidad del Pa\'{\i}s Vasco\\Aptdo. 644, 48080 Bilbao, Spain}
\email[J. C. Peral]{juancarlos.peral@ehu.es}
\thanks{A. D. was supported by the Ministry of Science, Education and Sports, Republic of Croatia,
grant 037-0372781-2821. J. C. P. was supported by the UPV/EHU grant EHU 10/05.}
\begin{abstract}
 There are $26$ possibilities for the  torsion group of elliptic curves defined over  quadratic number fields. We present examples of high rank elliptic curves with given torsion group which give the current records for most of the torsion groups.
In particular, we show that for each possible torsion group, except maybe for $\Z/15\Z$,
there exist an elliptic curve over some quadratic field with this torsion group and with rank $\geq 2$.
\end{abstract}

\maketitle

\section{Introduction }
Let $K=\Q(\sqrt{d})$, $d\in\Z$, be a quadratic number field and $E$ an elliptic curve defined over $K$. Then the torsion subgroup for $E$ is one of:
\begin{align*}
\Z/m\Z, &\quad 1\le m \le 18, m\neq 17,\\
\Z/2\Z \times\Z/2m \Z, &\quad1\le m\le 6,\\
\Z/3\Z \times\Z/3m \Z, &\quad m= 1,2, \quad \hbox{only possible for $\Q(\sqrt{-3})$,}\\
\Z/4\Z \times\Z/4\Z, & \quad \hbox{only possible for $\Q(\sqrt{-1})$}
\end{align*}
(see \cite{K} and  \cite{KM}).
In \cite{R}, Rabarison gave a parametrization of the curves having such torsion groups.
His results together with recent results from \cite{BBDN} show that for each possible
torsion group there exist an elliptic curve over some quadratic field with this torsion group and positive rank. In this paper, we present examples of elliptic curves having high rank and given torsion group over a quadratic field. For most of the torsion groups our results give the current records.
In particular, our results show that for each possible torsion group,
except maybe for $\Z/15\Z$, there exist an elliptic
curve over some quadratic field with this torsion group and with rank $\geq 2$.

For elliptic curves over $\Q$, by Mazur's theorem there are $15$ possible torsion groups.
The current rank records and the data on corresponding curves can be found in the table \cite{D}
(for previous records see e.g. \cite{KS}). It is known since 2003 that for each of these
$15$ torsion groups there exist a curve over $\Q$ with this torsion group and with rank $\geq 3$,
and this is still the best known result of that type. Namely, despite big efforts in last 10 years,
no elliptic curve over $\Q$ with torsion $\Z/2\Z \times \Z/8\Z$ and rank $\geq 4$ has been found so far.

\section{Current records} \label{sec:table}
We summarize the current records, including our results, in Table~\ref{t:table}.
For each torsion group $T$ from the above list, we define  $B(T)$ as the supremum
of the ranks of elliptic curves defined over any quadratic field and having torsion group $T$ (for the trivial torsion group we put $T=0$). The third column in Table~\ref{t:table} (l.b.) gives a lower bound of $B(T)$ for $T$ given in the first column, while the second column gives the discriminant $d$
of the corresponding quadratic number field $K$.
%For the comparison, we also give current records
%for curves over $\Q$, of course, only for the $15$ torsion groups which can appear over $\Q$
%by Mazur's theorem.
Numbers 4.x in the reference column correspond to the subsections of Section~\ref{examples} of the present paper, where the details of the corresponding curve appear.
The marks $*$ in the table refer to the curves obtained by an application of Mestre's general construction which produces quadratic fields with huge discriminants. In the case of torsion $\Z/2\Z\times \Z/2\Z$, the discriminant $d_2$ is too large to fit in the table, but will be given in Subsection \ref{z2z2}. Values in boldface in the third column are conditional, assuming the Parity Conjecture.

\begin{table}[htbp!]
\centering
\caption{Current records}\label{t:table}
\small
\begin{tabular}{@{} c c c c @{}}
\toprule
$T$        &       $d$ & l.b. & Reference\\
\midrule
$0$                    & $*$           & 30      & \cite{D,E,M1} \\
$\Z/2\Z$               & $-1$          & 28      & \ref{z2} \\
$\Z/3\Z$               & $ *$          & 15      & \cite{D,M1} \\
$\Z/4\Z$               & $-25689$      & 15      & \ref{z4} \\
$\Z/5\Z$               & $*$           & 10      & \cite{D,M1} \\
$\Z/6\Z$               & $3521$        & 11      & \ref{z6} \\
$\Z/7\Z$               & $*$           & 7       & \cite{D,M1} \\
$\Z/8\Z$               & $-227$        & 9       & \ref{z8} \\
$\Z/9\Z$               & $-155$        & 6       & \ref{z9} \\
$\Z/10\Z$              & $-2495$       & 7       & \ref{z10} \\
$\Z/11\Z$              & $-3239$       & 2       & \ref{z11} \\
$\Z/12\Z$              & $2014$        & 7       & \ref{z12} \\
$\Z/13\Z$              & $193$         & 2       & \cite{R} \\
$\Z/14\Z$              & $265$         & 2       & \ref{z14}, \cite{R} \\
$\Z/15\Z$              & $-7$          & 1       & \ref{z15}, \cite{BBDN} \\
$\Z/16\Z$              & $1785$        & 2       & \ref{z16} \\
$\Z/18\Z$              & $26521$       & 2       & \cite{BBDN} \\
$\Z/2\Z\times \Z/2\Z$  & $d_2$         & 19 \textbf{20} & \ref{z2z2} \\
$\Z/2\Z\times \Z/4\Z$  & $-83201$      & 13      & \ref{z2z4} \\
$\Z/2\Z\times \Z/6\Z$  & $624341$      & 10      & \ref{z2z6} \\
$\Z/2\Z\times \Z/8\Z$  & $31230597$    & 8       & \ref{z2z8} \\
$\Z/2\Z\times \Z/10\Z$ & $1065333545$  & 4 \textbf{5} & \ref{z2z10}, \cite{BBDN} \\
$\Z/2\Z\times \Z/12\Z$ & $2947271015$  & 4       & \cite{BBDN} \\
$\Z/3\Z\times \Z/3\Z$  & $-3$          & 7       & \cite{J} \\
$\Z/3\Z\times \Z/6\Z$  & $-3$          & 6       & \cite{J} \\
$\Z/4\Z\times \Z/4\Z$  & $-1$          & 7       & \cite{DJ} \\
\bottomrule
\end{tabular}
\end{table}

\section{Methods of search} \label{methods}
In the case of the $15$ possible torsion  groups of elliptic curves over $\Q$ (by Mazur's theorem), we consider curves with rational coefficients, and in order to determine their rank over a quadratic field $\Q(\sqrt{d})$ we use
the formula (\cite{B}, see also \cite{SZ} or \cite[Exercise 10.16]{S})
\begin{equation}\label{tworanks}
\rank(E(\Q(\sqrt{d}))=\rank(E(\Q))+\rank(E^{(d)}(\Q)),
\end{equation}
where $E^{(d)}$ denotes the $d$-quadratic twist of $E$.
One way to get high rank curves over quadratic fields with one of these $15$
torsion groups is to start with one of the record curves with
this torsion over $\Q$ (these curves are collected in the tables \cite{D}),
then search for a twist of this curve with high rank, and finally use \eqref{tworanks}.
However, the record curves usually have very large coefficients,
which makes it hard to compute the rank of its twist. Thus sometimes it is more profitable
to start with curves of rank slightly smaller than the corresponding record
(many such curves are also available in \cite{D}), and then use the full power
of standard sieving methods for finding high rank curves in parametric families of
elliptic curves. These methods include computations of Mestre-Nagao sums (they
are indicators of good candidates of high rank, see e.g. \cite{N}) and
the Selmer rank (it gives an upper bound for the rank which is sometimes sharp).
In some cases the root number is used either for reaching conditional results through the Parity Conjecture or as an additional sieving parameter.
We implemented these sieving procedures in \texttt{Pari} \cite{pari}. For computing the rank we used Cremona's program \texttt{mwrank} \cite{mwrank}.

For the curves with record rank over $\Q$, we were usually able to compute the rank of the twists 
only up to $d<10^4$ or $d<10^5$, so that gave the limitation for our search. 
On the other hand, for starting curves with small coefficients, we could search for high-rank 
twists up to $d<10^8$. In order to make this search practical, we have chosen strong conditions 
for the Mestre-Nagao sum. As an illustration, we can mention the $31230597$-twist of rank $5$ 
in the case of torsion group $\Z/2\Z\times \Z/8\Z$ (for details see Subsection \ref{z2z8}). 
Here the search procedure took around $10$ days of CPU time on an Intel Xeon E5430. 

We should mention a general result due to Mestre \cite{M1} (see also \cite{ST})
who proved that for any elliptic curve $E$ over $\Q$ there exist infinitely
 many quadratic twists with rank $\geq 2$.
For some torsion groups of odd order (trivial, $\Z/3\Z$, $\Z/5\Z$, $\Z/7\Z$) the largest known
ranks were obtained by applying Mestre's result \cite{M1} to the record curves over $\Q$.
Thus we will not give further details for these cases.
We mention also results announced by Mestre \cite{M2} and proved by Rubin and
Silverberg \cite{RS1,RS2} that if the torsion group of $E$ is
$\Z/2\Z\times \Z/4\Z$, resp. $\Z/2\Z\times \Z/8\Z$,
then there exist infinitely many quadratic twists of $E$ with rank $\geq 3$, resp. $\geq 4$.
However, our results for these two cases improve those obtained by direct application of
results from \cite{RS1,RS2} to the record curves over $\Q$.

For the torsion groups of the form $\Z/2\Z\times \Z/2m\Z$, we may take the record curve
with torsion $\Z/2m\Z$ over $\Q$, and then determine the quadratic field over which this curve
contains torsion $\Z/2\Z\times \Z/2m\Z$.

Formula \eqref{tworanks} can be also used for curves over a quadratic field, even if its torsion group is not one of the $15$ which are possible for curves over $\Q$, when they admit a model with rational coefficients. Indeed, this was applied in \cite{DJ,J}
for torsion groups $\Z/4\Z\times \Z/4\Z$, $\Z/3\Z\times \Z/3\Z$ and $\Z/3\Z\times \Z/6\Z$.
For the remaining torsion groups, we use the parametrizations given in \cite{R},
and try to compute the rank (or at least to get some information on it) for curves
corresponding to small parameters by using procedures available in \texttt{Magma} (\cite{Magma}).

Let us mention that for torsion groups $\Z/13\Z$ and $\Z/18\Z$, the current records with rank $2$  might be very hard to improve, since it was recently proved in \cite{BBDN} that all such curves  over quadratic fields have even rank.

\section{Examples} \label{examples}
%%%%
%\subsection{Curves with trivial torsion. $B(0, K)\ge 29$.} \label{z1}
%Elkies example of rank $28$ curve over Q with trivial torsion is
%birationally equivalent to
%\begin{align*}
%E_1:\quad y^2 =& x^3-\\
%&26007820090585882454203037839302751405205578484583927626187 x+\\
%&160877407990894564252589098013679906899400599794056833\\
%&8381502120424946520558796200800134.
%\end{align*}
%We have found a rank $1$  twist for the following value of  $d$:
%\begin{align*}
%d = &25793827509760603274338357335988598\\&15956399783437785323480000466351028118587.
%\end{align*}
%The equality \eqref{tworanks} gives the result.
%%%%%%%%
\subsection{Torsion $\Z/2\Z$. $B(\Z/2\Z)\ge 28$.} \label{z2}
For this torsion group the highest rank known over $\Q$ is an example of a curve with rank $19$ found by Elkies (see \cite{D}). For our purpose, it is more interesting to use an example due to Watkins (see \cite{ACP,D}), who found a curve with rank $14$ and $j$-invariant $1728$, namely
\[
y^2=x^3-402599774387690701016910427272483\,x.
\]
It is isomorphic to its $(-1)$-twist. Hence, the rank of the $(-1)$-twist
of this elliptic curve is also $14$, and now it follows from \eqref{tworanks} that
this curve has rank $28$ over $\Q(\sqrt{-1})$.
%%%%
%\subsection{Torsion $\Z_3$. $B(\Z_3, K)\ge 14$.}
%Eroshkin (see \cite{D}) found a curve with torsion group $\Z_3$ and rank $13$ over $\Q$. The curve is %given by
%\begin{align*}
%y^2 + x y = x^3-& 560715933702165990261993692150795879540 x
%+\\&5299428030171662962897867758309003693598430128674403539600.
%\end{align*}
%The twist for the following value of $d$
%\[
%d = 7574698155069228437966180387102961667517705
%\]
%has rank at least $1$ and conditioned by the Parity Conjecture $2$. The torsion is the same that over %$\Q$, so the result follows.
%%%%
\subsection{Torsion $\Z/4\Z$. $B(\Z/4\Z)\ge 15$.} \label{z4}
The curve with the highest known rank over $\Q$ having this torsion group is
an example with rank $12$ found by Elkies in 2006 (see \cite{D}).
Since its coefficients are very large, we have considered also curves with rank $9$, $10$ and $11$ from \cite{D}. In that way, we found several curves with rank $14$ over certain quadratic fields, and one example with rank $15$. This example was found starting with the rank $9$ curve:
\begin{multline*}
y^2 = x^3-x^2-63101137631999143241257265\,x \\
+ 191341468112384598938400129107933347137
\end{multline*}
found by Eroshkin in 2008.
Its twist by $d=-25689$ has rank $6$, so an application of \eqref{tworanks}
gives that $B(\Z/4\Z)\ge 15$.
%%%%
%  \subsection{Torsion $\Z_5$. $B(\Z_5, K)\ge 9$.}
%  Dujella and Lecaheux \cite{D} found the elliptic curve
%\begin{align*}
%  E_5:\quad  y^2 + x y = &x^3-1346404541224901580948752030 x+\\
%&19318553476047119421649468184366353981476
%\end{align*}
% with torsion group $\Z_5$ and rank $8$ over $\Q$.
%
% Its twist by $d=198908838922712275879496562795$ has rank $1$ and torsion group $\Z_5$ over %$\Q(\sqrt{198908838922712275879496562795})$ so the result follows.
% %%%%
\subsection{Torsion $\Z/6\Z$. $B(\Z/6\Z)\ge 11$.} \label{z6}
In this case the highest known rank over $\Q$  known for this torsion group is  $8$
(Dujella, Eroshkin, Elkies (2008), see \cite{D}).
We got the best result starting with a rank $7$ curve, found by Dujella in 2001, which is given by
$$
y^2 + x\,y = x^3 - x^2 - 45123702275641081919424\,x
+ 936989213947498862436000.
$$
Its $3521$-twist has rank at least $4$ so the result follows again  by~\eqref{tworanks}.
%%%%%
%   \subsection{Torsion $\Z_7$. $B(\Z_7, K)\ge 6$.}  Dujella and Kulesz found the elliptic curve
%   \begin{align*}
%  E_7:\quad  y^2 + x y =& x^3-3871874529581532682991685 x
%+\\&2932443695831973034287157248711276225
% \end{align*}
% with rank $5$ and torsion group $\Z_7$. Its twist by $d=263$ has rank at least $1$ and the same torsion %group so we have the result.
% %%%%%
\subsection{Torsion $\Z/8\Z$. $B(\Z/8\Z)\ge 9$.} \label{z8}
For this torsion group the curve with highest known rank was found by Elkies in 2006 (see \cite{D}) and has rank~$6$. For our search we use the curve with rank $5$, found by Dujella \& Lecacheux in 2002, which is given by
%{\small
\begin{multline*}
y^2= x^3 + x^2 - 11849634571550798667743047864720\,x  \\
+ 15613761915399875450490670165233536220551598068.
\end{multline*}
%}%
Since its twist by $d=-227$ has rank $4$, the result follows.
%%%%%
\subsection{Torsion $\Z/9\Z$. $B(\Z/9\Z)\ge 6$.} \label{z9}
Fisher found in 2009 a curve with torsion $\Z/9\Z$ and rank $4$ over $\Q$ (see \cite{D}).
By applying the general result by Mestre, we get a curve of rank $6$ over a quadratic field. However, this field has a huge discriminant, and an example of rank~$6$ with much smaller discriminant can be obtained starting with the elliptic curve
$$
y^2 + x\,y + y = x^3- x^2 -41368267697099\,x
+ 102411668493915101147,
$$
of rank $3$ over $\Q$ (Dujella 2001), and noting that its twist by $d=-155$ also has rank $3$ over $\Q$. Thus it has rank $6$ over $\Q(\sqrt{-155})$ by~\eqref{tworanks}.
%%%%%

\subsection{Torsion $\Z/10\Z$. $B(\Z/10\Z)\ge 7$.}  \label{z10}
We use the curve with rank $4$ and torsion $\Z/10\Z$ which is the rank record for this torsion group (Dujella 2004, see \cite{D}). The  curve is
%{\small
\begin{multline*}
y^2 + x\,y = x^3 - 127381738643041574974581021420318985\,x \\
+ 17495594046612039766866496413577998621609407092547225.
\end{multline*}
%}%
Its twist by $d=-2495$ has rank  equal to $3$, so the result follows from~\eqref{tworanks}.

Another examples with rank $7$ and torsion $\Z/10\Z$ over a quadratic field can be
obtained by starting with curves with smaller rank. E.g.
the curve
$$
y^2 + x\,y = x^3 -5313234280\,x  + 149068288642400
$$
(Womack 2000, see \cite{W})
has torsion $\Z/10\Z$, rank $2$ over $\Q$ and its $318855485$-twist has rank $5$, while the curve
\begin{multline*}
y^2 + x\,y = x^3 - 4281263352573652971565\,x \\
+ 107821663654697042219512111579217
\end{multline*}
(Dujella 2001, see \cite{D})
has torsion $\Z/10\Z$, rank $3$ over $\Q$ and its $3007$-twist has rank $4$.
%%%%%%

\subsection{Torsion $\Z/11\Z$. $B(\Z/11\Z)\ge 2$.} \label{z11}
This torsion cannot appear over $\Q$. Rabarison \cite{R} showed the existence of curves with torsion group $\Z/11\Z$ and rank at least $1$ over $\Q(\sqrt{d})$ for $d= -3239,-599,-47,6,7,22,73$ and $193$.
The curve over $\Q(\sqrt{-3239})$ is
$$
y^2 + \frac{5\sqrt{-3239} - 796}{9}\,x\,y + (45\sqrt{-3239} + 855)y
= x^3 + (45\sqrt{-3239} + 855)x^2.
$$
Rabarison stated that its rank over $\Q(\sqrt{-3239})$ is $1$.
Since it has a point of infinite order ($[1-\sqrt{-3239}, 2\sqrt{-3239}-1082]$)
and root number equal to $1$, the Parity Conjecture implies that its rank should be at least~$2$. Unfortunately, we were not able to find the second independent point.
We have also an example over $\Q(\sqrt{1129})$ with conditional rank $\geq 2$.
However, we were able to find one example with torsion $\Z/11\Z$ and rank $\geq 2$
over $\Q(\sqrt{561})$, unconditionally.
Namely, the curve
$$
y^2 + \frac{-10\sqrt{561} + 893}{1008}\,x\,y + \frac{-35\sqrt{561} + 210}{10368}y
= x^3 + \frac{-35\sqrt{561} + 210}{10368} x^2
$$
has the torsion group $\Z/11\Z$ over $\Q(\sqrt{561})$
and the following two independent points of infinite order:
$\left[ \frac{-\sqrt{561} - 869}{6912}, \frac{1709\sqrt{561} + 147421}{6967296} \right]$ and
$\left[ \frac{5\sqrt{561} - 65}{1344}, \frac{-1285\sqrt{561} + 32385}{2032128} \right]$.
%%%%%%%

\subsection{Torsion $\Z/12\Z$. $B(\Z/12\Z)\ge 7$.}  \label{z12}
Although the highest known rank for this torsion is $4$ (Fisher 2008, see \cite{D}),
we apply our procedure to several curves with rank $2$ and $3$.
For the curve with rank $3$ and torsion $\Z/12\Z$ found by Rathbun in 2003 (see \cite{D}) given by
%{\small
\begin{multline*}
y^2 + x\,y + y = x^3- x^2 -42403753582533569425032932\,x \\
+ 106274144228004532427905140464314177031,
\end{multline*}
%}%
we find that its twist by $d=-106071$ has rank $4$. Another curve with rank $3$ and torsion
$\Z/12\Z$ is
%{\small
\begin{multline*}
y^2 + x\,y = x^3 - 544753256053055692212823356675\,x \\
+ 154756127532691562955214620687209364995464257,
\end{multline*}
%}%
found by Dujella in 2005. We find that its $2905$-twist has rank $4$. Finally, we consider a curve
with rank $2$:
\begin{multline*}
y^2 + x\,y = x^3-x^2-136659485377389900024\,x \\
+ 612767297917647098548240331268,
\end{multline*}
which has $2014$-twist with rank $5$. Again, an application of \eqref{tworanks} gives the announced result.
%%%%%%%
%       \subsection{Torsion $\Z_{13}$. $B(\Z_{13}, K)\ge 2$.}  In Rabarison (theorem $6$, \cite{R}), it is %proved that for the field $\Q(\sqrt{193})$ there exist an elliptic curve with torsion group $\Z_{13}$  %and rank $2$ over this field.
%%%%%%%%%
\subsection{Torsion $\Z/14\Z$. $B(\Z/14\Z)\ge 2$.}  \label{z14}
In \cite{R}, Rabarison gave examples of elliptic curves with torsion group $\Z/14\Z$
 and rank $\geq 1$ over the quadratic fields $\Q(\sqrt{-23})$ and $\Q(\sqrt{265})$.
The curve over $\Q(\sqrt{265})$ is
$$
y^2 + \frac{-18\sqrt{265} + 283}{145}\,x\,y
+ \frac{726\sqrt{265} - 12990}{21025}\,y
= x^3 + \frac{726\sqrt{265} - 12990}{21025}\,x^2.
$$
Using \texttt{Magma} \cite{Magma}, we find that this curve has rank equal to $2$, with independent points with $x$-coordinates
\[
\frac{-133118\sqrt{265} + 4584710}{36395145} \,\,\text{ and }\,\,
\frac{-246\sqrt{265} + 4902}{4205}.
\]
%%%%%%%

\subsection{Torsion $\Z/15\Z$. $B(\Z/15\Z)\ge 1$.} \label{z15}
In \cite{BBDN}, an example with torsion group $\Z/15\Z$ and positive rank is given over
$\Q(\sqrt{345})$. Here we give such an example over $\Q(\sqrt{-7})$. This is the curve
$$
y^2 + (-2\sqrt{-7}+15)x\,y + (26\sqrt{-7} -14)y
= x^3 + (26\sqrt{-7} -14)x^2,
$$
obtained by inserting
\[t=(-\sqrt{-7}-5)/8,\quad s=(\sqrt{-7} -11)/16\]
in the general parametrization of curves with torsion $\Z/15\Z$ given in \cite{R}.
This curve has a point of infinite order $[6\sqrt{-7} - 98, 136\sqrt{-7} + 1064]$
over $\Q(\sqrt{-7})$.
%%%%%%%
\subsection{Torsion $\Z/16\Z$. $B(\Z/16\Z)\ge 2$.} \label{z16}
In \cite{R}, Rabarison gave an example of elliptic curve with torsion $\Z/16\Z$ and
rank $\geq 1$ over $\Q(\sqrt{10})$. We improve this result by finding an example with
torsion $\Z/16\Z$ and rank equal to $2$ over $\Q(\sqrt{1785})$. The curve is
{\small
$$
y^2 + \frac{-247\sqrt{1785} + 80605}{71680}\,x\,y
+ \frac{-203\sqrt{1785} - 11271}{917504}\,y
= x^3 + \frac{-203\sqrt{1785} - 11271}{917504}\,x^2.
$$
}%
Using \texttt{Magma} \cite{Magma}, we find that this curve has rank~$2$ over $\Q(\sqrt{1785})$, with independent points with $x$-coordinates
\[\frac{10763581607\sqrt{1785} + 456934229139}{894936283136} \,\,\text{ and }\,\,
\frac{11\sqrt{1785} + 663}{50540}.\]
%%%%%
%          \subsection{Torsion $\Z_{18}$. $B(\Z_{18}, K)\ge 2$.} This result is proved in \cite{BBDN}.
%
% %%%%%%%%%%%%%%%%%%%

\subsection{Torsion $\Z/2\Z \times \Z/2\Z$. $B(\Z/2\Z \times \Z/2\Z)\ge 19$.}\label{z2z2}
 In this case the argument is different: we start with  a curve having torsion
 $\Z/2\Z$ over $\Q$ and rank $19$. This curve has the record rank for this torsion group and was found by Elkies in 2009 (see  \cite{D}).
 We write it in the form
\[y^2=(x+655364911965298267755181)g(x),\]
where $g$ is a monic polynomial of degree~$2$ with integer coefficients whose discriminant is
\[d_2=-3901785498412536920668361993073821511.\]
The curve has two additional $2$-torsion points over $\Q(\sqrt{d_2})$ and its torsion group over this field is  $\Z/2\Z \times \Z/2\Z$. The rank is at least $19$ over $\Q(\sqrt{d_2})$. Since the root number for the $d_2$-twist of the curve is $-1$, assuming the Parity Conjecture, we have that the rank should be at least $20$.

%%%%%%%%%%%%%%%%
\subsection{Torsion $\Z/2\Z \times \Z/4\Z$. $B(\Z/2\Z \times \Z/4\Z)\ge 13$.}  \label{z2z4}
The record rank over $\Q$ for this torsion is $8$, due to
Elkies (2005), Eroshkin (2008), Dujella \& Eroshkin (2008) (see \cite{D}).
One of the known rank-record curves, found by Eroshkin, is given by
%{\small
\begin{multline*}
y^2 = x^3 + x^2-23686061832482481624168232900\,x \\
+ 1401294826072670363740983663536729053022048.
\end{multline*}
%}%
Its twists by $d=-83201$ and $d=109499$ both have rank $5$, so the result follows.

Another example is given by the curve
%{\small
\begin{multline*}
y^2 = x^3 - x^2 - 866893152450363503763740085700\,x \\
+61220734062068506723288644020689511073795652,
\end{multline*}
%}%
found by Dujella \& Eroshkin. It has rank $8$, while its $117589$-twist has rank $5$.

%%%%%%%%%%%%%%%%
\subsection{Torsion $\Z/2\Z \times \Z/6\Z$. $B(\Z/2\Z \times \Z/6\Z)\ge 10$.} \label{z2z6}
Although the record rank over $\Q$ for this torsion is $6$ (Elkies, 2006),
we use a curve having this torsion group and rank $5$, found by
Dujella \& Lecacheux in 2002 (see \cite{D}). The curve is given by
%{\small
\begin{multline*}
y^2 + x\,y = x^3 - 2353799432200918732090882185\,x \\
+ 39989567111692230080439457690563434811225.
\end{multline*}
%}%
The result follows from the fact that its $624341$-twist has rank~$5$ over $\Q$.

%%%%%%%%%%%%%
\subsection{Torsion $\Z/2\Z \times \Z/8\Z$. $B(\Z/2\Z \times \Z/8\Z)\ge 8$.} \label{z2z8}
The highest known rank over $\Q$ for curves with this torsion is $3$.
In total, there are $27$ such curves known, found by many authors (see \cite{D}).
We use the curve with rank $3$, found independently by Connell and Dujella in 2000:
$$
y^2 + x\,y = x^3 - 15745932530829089880\,x
+ 24028219957095969426339278400.
$$
The result is a consequence of the fact that its twist by $d=31230597$ has rank~$5$.
%%%%%%%%%%%%%%%%%
\subsection{Torsion $\Z/2\Z \times \Z/10\Z$. $B(\Z/2\Z \times \Z/10\Z)\ge 4$.} \label{z2z10}
This torsion group does not appear over $\Q$, so we apply an argument similar to the one used in Subsection \ref{z2z2}. We start from a curve with rank $4$ over $\Q$ and torsion $\Z/10\Z$ and we look for  twists with additional torsion points. The starting curve is
%{\small
\begin{multline*}
y^2 + x\,y = x^3- 186734461375182851611482374885236743900\,x\\
+ 979124235576847382684597033191041431087287234552438410000,
\end{multline*}
%}%
found by Dujella in 2008 (see \cite{D}).
We find that it has torsion $\Z/2\Z \times \Z/10\Z$ over
the quadratic field $\Q(\sqrt {1065333545} )$, so the result follows.
Let us mention that  a curve of rank $\geq 4$ and torsion $\Z/2\Z \times \Z/10\Z$ over $\Q(\sqrt {55325286553})$ was found in \cite{BBDN} by the same method. However, here the corresponding twist has the root number $-1$, so assuming the Parity Conjecture, this curve has rank $\geq 5$.

%%%%%%%%%%%%%%%%%
%\subsection{Torsion $\Z_2 \times \Z_{12}$. $B(\Z_2 \times \Z_{12}, K)\ge 4$.} \label{z2z12}
% This result is proved in \cite{BBDN}.
%%%%%%%%%%%%%%%%%
%\subsection{Torsion $\Z_3 \times \Z_{3}$. $B(\Z_3 \times \Z_3, K)\ge 7$.}  The torsion considered in this %subsection only appears over $\Q(\sqrt{-3})$ and  the  result it is proved in \cite{J}.
%%%%%%%%%%%
%\subsection{Torsion $\Z_3 \times \Z_6$. $B(\Z_3 \times \Z_3, K)\ge 6$.}  The torsion considered in this %subsection only appears over $\Q(\sqrt{-3})$ and  the  result it is proved in \cite{J}.
%%%%%%%%%%%
%\subsection{Torsion $\Z_4 \times \Z_4$. $B(\Z_4 \times \Z_4, K)\ge 7$.}  The torsion considered in this %subsection only appears over $\Q(\sqrt{-1})$ and  the  result it is proved in \cite{DJ}.

\bigskip

{\bf Acknowledgement}. The authors would like to thank Filip Najman and the referee for very useful
comments on the previous version of this paper.

 \end{document}